\documentclass[10pt]{amsart}
\usepackage{epsfig,url}
\usepackage{mathdots}

\newtheorem{theorem}{Theorem}[section]

\theoremstyle{definition}

\newtheorem{note}[theorem]{Note}

\theoremstyle{remark}

\begin{document}

\title[Pochhammer symbol with negative indices]{Pochhammer symbol with negative indices.  \\ A new rule for the method of brackets. }

\author{Ivan Gonzalez}
\address{Instituto  de Fisica y Astronomia \\ Universidad de Valparaiso \\ Chile}
\email{ivan.gonzalez@uv.cl}

\author{Lin Jiu}
\address{Department of Mathematics,
Tulane University, New Orleans, LA 70118}
\email{ljiu@tulane.edu}

\author{Victor H. Moll}
\address{Department of Mathematics,
Tulane University, New Orleans, LA 70118}
\email{vhm@tulane.edu}

\subjclass[2010]{Primary 33C99, Secondary 33C10}

\date{\today}

\keywords{ Method of brackets, Pochhammer symbol, Bessel integral}

\begin{abstract}
The method of brackets is a method of integration based upon  a small number of 
heuristic rules.  Some of these have been made rigorous. An example of an integral involving the 
Bessel function is used to motivate a  new evaluation rule.
\end{abstract}

\maketitle

\newcommand{\ba}{\begin{eqnarray}}
\newcommand{\ea}{\end{eqnarray}}
\newcommand{\ift}{\int_{0}^{\infty}}
\newcommand{\nn}{\nonumber}
\newcommand{\no}{\noindent}
\newcommand{\lf}{\left\lfloor}
\newcommand{\rf}{\right\rfloor}
\newcommand{\realpart}{\mathop{\rm Re}\nolimits}
\newcommand{\imagpart}{\mathop{\rm Im}\nolimits}

\newcommand{\op}[1]{\ensuremath{\operatorname{#1}}}
\newcommand{\pFq}[5]{\ensuremath{{}_{#1}F_{#2} \left( \genfrac{}{}{0pt}{}{#3}
{#4} \bigg| {#5} \right)}}

\newtheorem{Definition}{\bf Definition}[section]
\newtheorem{Thm}[Definition]{\bf Theorem}
\newtheorem{Example}[Definition]{\bf Example}
\newtheorem{Lem}[Definition]{\bf Lemma}
\newtheorem{Cor}[Definition]{\bf Corollary}
\newtheorem{Prop}[Definition]{\bf Proposition}
\numberwithin{equation}{section}

\section{Introduction}
\label{sec-intro}

The evaluation of definite integrals is connected, in a surprising manner, to many topics in Mathematics. The last
author has described in \cite{moll-2002a} and \cite{moll-2010c} some of these connections.   Many of these evaluations 
appear in the classical table \cite{gradshteyn-2015a} and  proofs of these entries have appeared in a 
series of papers starting with \cite{moll-2007a} and the latest one is \cite{amdeberhan-2015c}.  An interesting 
new method of integration, developed in \cite{gonzalez-2007a} in the context of integrals coming from Feynman diagrams, 
is illustrated here in the evaluation of entry $6.671.7$ in \cite{gradshteyn-2015a}:
\begin{equation}
I:=\int_{0}^{\infty} J_{0}(ax) \sin(bx) \, dx = \begin{cases} 
0, & \quad \text{ if } 0< b < a, \\
1/\sqrt{b^{2}-a^{2}}, & \quad \text{ if } 0 < a < b.
\end{cases}
\label{entry1}
\end{equation}

This so-called \textit{method of brackets} is based upon a small number of heuristic rules. These 
are described in the next section.  Sections 
\ref{gr-entry1} and \ref{sec-alternative} present the evaluation of \eqref{entry1}. The conflict between these two 
evaluations is resolved in Section \ref{sec-extension} with the proposal of a new rule dealing with the extension 
of the Pochhammer symbol to negative integer indices.

\section{The method of brackets}
\label{sec-method}

 A method to evaluate integrals 
over the half line $[0, \, \infty)$, based on a small number of rules has 
been developed in \cite{gonzalez-2007a,gonzalez-2008a,gonzalez-2009a}. This 
method of brackets is described 
next. The heuristic rules are currently being made rigorous in
 \cite{amdeberhan-2012b} and \cite{jiu-2015e}. The reader will find in 
\cite{gonzalez-2014a,gonzalez-2010a,gonzalez-2010b} a large collection of 
evaluations of definite integrals that illustrate the power and flexibility 
of this method. 

For $a \in \mathbb{C}$, the symbol 
\begin{equation}
\langle a \rangle  \mapsto \ift x^{a-1} \, dx 
\end{equation}
is the {\em bracket} associated to the (divergent) integral on the right. The 
symbol 
\begin{equation}
\phi_{n} := \frac{(-1)^{n}}{\Gamma(n+1)}
\end{equation}
\noindent
is called the {\em indicator} associated to the index $n$. The notation
$\phi_{i_{1}i_{2}\cdots i_{r}}$, or simply 
$\phi_{12 \cdots r}$, denotes the product 
$\phi_{i_{1}} \phi_{i_{2}} \cdots \phi_{i_{r}}$. 

\medskip

\noindent
{\bf {\em Rules for the production of bracket series}}

\smallskip

\noindent
${\mathbf{Rule \, \, P_{1}}}$.  If the function $f$ is given by the formal
power series
\begin{equation}
f(x)=\sum_{n=0}^{\infty}a_{n}x^{\alpha n+\beta-1},
\end{equation}
then the improper integral of $f$ over the positive real line is
formally written as the \emph{bracket series }
\begin{equation}
\int_{0}^{\infty}f(x)dx=\sum_{n}a_{n}\left\langle \alpha n+\beta\right\rangle .
\end{equation}

\smallskip

\noindent
${\mathbf{Rule \, \, P_{2}}}$. 
For $\alpha \in \mathbb{C}$, the multinomial power 
$(a_{1} + a_{2} + \cdots + a_{r})^{\alpha}$ is assigned the 
$r$-dimension bracket series 
\begin{equation}
\sum_{n_{1}} \sum_{n_{2}}  \cdots \sum_{n_{r}}
\phi_{n_{1}\, n_{2} \,  \cdots n_{r}}
a_{1}^{n_{1}} \cdots a_{r}^{n_{r}} 
\frac{\langle -\alpha + n_{1} + \cdots + n_{r} \rangle}{\Gamma(-\alpha)}.
\end{equation}

\noindent
{\bf {\em Rules for the evaluation of a bracket series}}

\smallskip

\noindent
${\mathbf{Rule \, \, E_{1}}}$. 
The one-dimensional bracket series is assigned the value 
\begin{equation}
\sum_{n} \phi_{n} f(n) \langle an + b \rangle = 
\frac{1}{|a|} f(n^{*}) \Gamma(-n^{*}),
\end{equation}
\noindent
where $n^{*}$ is obtained from the vanishing of the bracket; that is, $n^{*}$ 
solves $an+b = 0$. 

\smallskip

The next rule provides a value for multi-dimensional bracket series where 
 the number of sums is equal to the number of brackets. 

\smallskip

\noindent
${\mathbf{Rule \, \, E_{2}}}$. 
Assume the matrix
$A = (a_{ij})$ is non-singular, then the assignment is 
\begin{equation}
\sum_{n_{1}} \cdots \sum_{n_{r}} \phi_{n_{1} \cdots n_{r}} 
f(n_{1},\cdots,n_{r}) 
\langle a_{11}n_{1} + \cdots + a_{1r}n_{r} + c_{1} \rangle \cdots
\langle a_{r1}n_{1} + \cdots + a_{rr}n_{r} + c_{r} \rangle 
\nonumber
\end{equation}
\begin{equation}
=  \frac{1}{| \text{det}(A) |} f(n_{1}^{*}, \cdots n_{r}^{*}) 
\Gamma(-n_{1}^{*}) \cdots \Gamma(-n_{r}^{*}) 
\nonumber
\end{equation}
\noindent
where $\{ n_{i}^{*} \}$ 
is the (unique) solution of the linear system obtained from the vanishing of 
the brackets. There is no assignment if $A$ is singular. 

\smallskip

\noindent
${\mathbf{Rule \, \, E_{3}}}$. 
Each representation of an integral by a bracket series has
associated an {\em index of the representation} via 
\begin{equation}
\text{index } = \text{number of sums } - \text{ number of brackets}.
\end{equation}
\noindent 
It is important to observe that the index is attached to a specific 
representation of the integral and not just to integral itself.  The 
experience obtained by the authors using this method suggests that, among 
all representations of an integral as a bracket series, the one with 
{\em minimal index} should be chosen. 

The value of a multi-dimensional bracket series of positive index
is obtained by computing all the contributions of maximal rank 
by Rule $E_{2}$. These
contributions to the integral appear as series in the 
free parameters. Series converging in a
common region are added and divergent series are discarded. 

\smallskip

\begin{note}
A systematic procedure in the simplification of the series obtained by this procedure  has been used throughout the 
literature: express factorials in terms of the gamma function and the transform quotients of gamma terms into Pochhammer symbol,  defined by
\begin{equation}
(a)_{k} = \frac{\Gamma(a+k)}{\Gamma(a)}.
\end{equation}
\noindent
Any presence of a Pochhammer with a negative index $k$ is transformed by the rule 
\begin{equation}
(a)_{-k} = \frac{(-1)^{k}}{(1-a)_{k}}.
\label{rule-11}
\end{equation}
\end{note}

 The example discussed in the next two 
section provides motivation for an additional evaluation rule for the method of brackets.

\section{A first  evaluation of entry $\mathbf{6.671.7}$ in  Gradshteyn and Ryzhik}
\label{gr-entry1}

The evaluation of \eqref{entry1} uses the series 
\begin{equation}
J_{0}(ax) = \sum_{m=0}^{\infty}  \phi_{m} \frac{a^{2m}}{\Gamma(m+1) 2^{2m}} x^{2m}
\end{equation}
\noindent
and 
\begin{equation}
\sin(bx) = \sum_{n=0}^{\infty}  \phi_{n} \frac{\Gamma(n+1)}{\Gamma( 2 n + 2)} b^{2 n +1} x^{2 n+1}.
\end{equation}
\noindent
Therefore the integral in \eqref{entry1} is given by 
\begin{equation}
I = \sum_{m,n} \phi_{m,n} \frac{a^{2m} b^{2n+1} \Gamma(n+1)}{2^{2m} \Gamma(m+1) \Gamma(2 n + 2)}
\langle 2 m + 2n + 2 \rangle.
\end{equation}
\noindent 
The duplication formula for the gamma function transforms this expression to 
\begin{equation}
I = \frac{\sqrt{\pi}}{2} \sum_{m,n} \phi_{m,n} \frac{a^{2m} b^{2n+1} }{2^{2m+2 n} \Gamma(m+1) \Gamma( n+ 3/2)}
\langle 2 m + 2n + 2 \rangle.
\end{equation}

Eliminating the parameter $n$ using Rule $E_{1}$ gives $n^{*} = - m -1$ and produces 
\begin{eqnarray*}
I & = &  \frac{\sqrt{\pi}}{b} \sum_{m=0}^{\infty}  \phi_{m} \frac{1}{\Gamma \left(-m + \tfrac{1}{2} \right)} \left( \frac{a}{b} \right)^{2m} \\
& = & \frac{1}{b} \sum_{m = 0}^{\infty} \frac{(-1)^{m}}{m!}  \frac{ \left( \frac{a^{2}}{b^{2}} \right)^{m}}{\left( \tfrac{1}{2} \right)_{-m} }\\
& = & \frac{1}{b} \sum_{m=0}^{\infty} \left( \frac{1}{2} \right)_{m} \frac{1}{m!} \left( \frac{a^{2}}{b^{2}} \right)^{m} \\
& = & \frac{1}{\sqrt{b^{2} - a^{2}}}.
\end{eqnarray*}
\noindent 
The condition $|b|>|a|$ is imposed to guarantee the convergence of the series on the third line of the previous argument. 

\smallskip

The series obtained by eliminating the parameter $m$ by $m^{*} = -n -1$ vanishes because of the factor $\Gamma(m+1)$ in 
the denominator. The formula has been established.

\section{An alternative evaluation}
\label{sec-alternative}

A second evaluation of \eqref{entry1}  begins with 
\begin{equation}
\int_{0}^{\infty} J_{0}(ax) \sin(bx) \, dx = \sum_{m,n} \phi_{m,n} \frac{a^{2m} b^{2n+1}}{2^{2m} m! (2n+1)!} 
\langle 2m+2n+2 \rangle.
\end{equation}
\noindent
The evaluation of the bracket series is described next. 

\smallskip

\noindent
\textit{Case 1}. Choose $n$ as the free parameter. Then $m^{*} = -n-1$ and the contribution to the integral is 
\begin{equation}
I_{1} = \frac{1}{2} \sum_{n=0}^{\infty} \phi_{n} 
\left( \frac{a}{2} \right)^{-2n-2} b^{2n+1} \frac{1}{\Gamma(-n)} \frac{n!}{(2n+1)!} \Gamma(n+1).
\end{equation}
\noindent
Each term in the sum vanishes because the gamma function has a pole at the negative integers. 

\noindent
\textit{Case 2}. Choose $m$ as a free parameter. Then $2m+2n+2=0$ gives $n^{*} = -m-1$. The contribution to the 
integral is 
\begin{equation}
I_{2} = \frac{1}{2} \sum_{m} \phi_{m} \left( \frac{a}{2} \right)^{2m} b^{-2m-1} \frac{1}{m!} \cdot 
\frac{n!}{(2n+1)!}\Big{|}_{n=-m-1} \Gamma(m+1).
\label{int-1}
\end{equation}
\noindent
Now write 
\begin{equation}
\frac{n!}{(2n+1)!} = \frac{1}{(n+1)_{n+1}}
\end{equation}
\noindent
and \eqref{int-1} becomes 
\begin{equation}
I_{2} = \frac{1}{2b} \sum_{m=0}^{\infty} \frac{(-1)^{m}}{m!} \left( \frac{a}{2b} \right)^{2m} \frac{1}{(-m)_{-m}}.
\label{int-2}
\end{equation}
\noindent
Transforming the term $(-m)_{-m}$ by using 
\begin{equation}
(x)_{-n} = \frac{(-1)^{n}}{(1-x)_{n}}  \label{ivan-rule1}
\end{equation}
 gives 
\begin{equation}
(-m)_{-m}  = \frac{(-1)^{m}}{(1+m)_{m}}.
\end{equation}
\noindent
Replacing in \eqref{int-2} produces 
\begin{eqnarray*}
I_{2}  & = &   \frac{1}{2b} \sum_{k=0}^{\infty} \left( \frac{1}{2} \right)_{k} \frac{1}{k!} \left( \frac{a^{2}}{b^{2}} \right)^{k}  \\
 & = & \frac{1}{2} \frac{1}{\sqrt{b^{2}-a^{2}}}.
 \end{eqnarray*}

The method of brackets produces \textit{half of the expected answer}. Naturally, it is possible that the entry in 
\cite{gradshteyn-2015a} is erroneous (this happens once in a while). Some  numerical  computations  and the 
evaluation in the previous section, should convince
 the reader that this is not the case. The source of the error is the use of \eqref{ivan-rule1} for the 
  evaluation of the term $(-m)_{-m}$.  A discussion is presented in Section \ref{sec-extension}.

\section{Extensions of the Pochhammer symbol}
\label{sec-extension}

Rule $E_{1}$ of the  method of brackets requires  the evaluation of $f(n^{*})$. In many instances, this involves the 
evaluation of the Pochhammer symbol $(x)_{m}$ for $m \not \in \mathbb{N}$. In particular, the question of the value 
\begin{equation}
(-m)_{-m} \quad \text{ for }  \quad m \in \mathbb{N}
\end{equation}
\noindent
is at the core of the missing factor of $2$ in Section \ref{sec-alternative}.

The first extension of $(x)_{m}$ to negative values of $n$ comes from the identity 
\begin{equation}
(x)_{-m} = \frac{(-1)^{m}}{(1-x)_{m}}.
\label{poch-neg}
\end{equation}
\noindent
This is obtained from 
\begin{equation}
(x)_{-m}  = \frac{\Gamma(x-m)}{\Gamma(x)} = \frac{\Gamma(x-m)}{(x-1)(x-2) \cdots (x-m) \Gamma(x-m)}
\end{equation}
\noindent
and then changing the signs of each of the factors.  This is valid as long as $x$ is not a negative integer. The limiting 
value of the right-hand side in \eqref{poch-neg} as $x \to -km$, with $k \in \mathbb{N}$, is
\begin{equation}
(-km)_{-m} = \frac{(-1)^{m} \, (km)!}{((k+1)m)!}.
\end{equation}

On the other hand, the limiting value of the left-hand side is 
\begin{eqnarray*}
\lim\limits_{\varepsilon \to 0} \left( -k(m+ \varepsilon) \right)_{-(m+ \varepsilon)} & = & 
\lim\limits_{\varepsilon \to 0} \frac{\Gamma(-(k+1)m -(k+1) \varepsilon)}{\Gamma(-km - k  \varepsilon )} \\
& = & \lim\limits_{\varepsilon \to 0} 
\frac{\Gamma( -(k+1) \varepsilon) (-(k+1) \varepsilon)_{-(k+1)m} }
{\Gamma(- k \varepsilon) ( - k \varepsilon)_{-km}} \\
& = &  \lim\limits_{\varepsilon \to 0} 
\frac{\Gamma( -(k+1) \varepsilon)}{\Gamma(-k \varepsilon)}  
\frac{(-1)^{(k+1)m}}{(1+ (k+1) \varepsilon)_{(k+1)m}} \cdot \frac{(1+ k \varepsilon)_{km}}{(-1)^{km} } \\
& = & \frac{(-1)^{m} (km)!}{((k+1)m)!} \cdot \frac{k}{k+1}.
\end{eqnarray*}

Therefore the function $(x)_{-m}$ is discontinuous at $x = -km$, with 
\begin{equation}
\frac{\text{Direct } (-km)_{-m}}{\text{Limiting } (-km)_{-m}} = \frac{k+1}{k}.
\end{equation}
\noindent
For $k=1$, this ratio becomes $2$. This explains the missing $\tfrac{1}{2}$ in the calculation in Section 
\ref{sec-alternative}. Therefore it is the discontinuity of \eqref{rule-11} at negative integer values of the 
variables, what is responsible for the error in the evaluation of the integral \eqref{entry1}.

\medskip

This example suggest that the rules of the method of brackets should be supplemented with an additional one: 

\smallskip

\noindent
${\mathbf{Rule \, \, E_{4}}}$.  Let $k \in \mathbb{N}$ be fixed. In the evaluation of series, the rule 
\begin{equation}
(-km)_{-m} =  \frac{k}{k+1} \, \frac{(-1)^{m} \, (km)!}{((k+1)m)!} 
\end{equation}
\noindent
must be used to eliminate Pochhammer symbols with negative index and negative integer base.

\smallskip

A variety of other examples confirm that this heuristic rule leads to correct evaluations.

\medskip

\noindent
\textbf{Acknowledgments}.  The last author acknowledges the partial 
support of NSF-DMS 1112656.  The second  author is a  graduate student,  partially supported by the same grant.


\end{document}